\documentclass[final,1p,times]{elsarticle}

\usepackage{amssymb}
 \usepackage{amsthm}
\usepackage{amscd}
\usepackage{amsmath}
\usepackage{amsfonts}
\usepackage{amssymb}
\usepackage{graphicx}
\usepackage{color}
\newtheorem{theorem}{Theorem}

\newtheorem{lemma}[theorem]{Lemma}

\usepackage{mathrsfs}
\usepackage{titletoc}

\newcommand{\D}{\Delta}
\newcommand{\ra}{\rightarrow}
\newcommand{\p}{\partial}
\newcommand{\f}{\frac}

\renewcommand{\f}{\frac}

\newcommand{\be}{\begin{equation}}
\newcommand{\ee}{\end{equation}}
\newcommand{\bea}{\begin{eqnarray}}
\newcommand{\eea}{\end{eqnarray}}
\newcommand{\bna}{\begin{eqnarray*}}
\newcommand{\ena}{\end{eqnarray*}}

\renewcommand{\le}{\left}
\newcommand{\ri}{\right}

\journal{***}

\begin{document}

\begin{frontmatter}

\title{A Trudinger-Moser inequality for conical metric in the unit ball}

\author{Yunyan Yang}
 \ead{yunyanyang@ruc.edu.cn}
\author{Xiaobao Zhu}
 \ead{zhuxiaobao@ruc.edu.cn}

 \address{School of Mathematics,
Renmin University of China, Beijing 100872, P. R. China}

\begin{abstract}
In this note, we prove a Trudinger-Moser inequality for conical metric in the unit ball. Precisely, let $\mathbb{B}$ be
the unit ball in $\mathbb{R}^N$ $(N\geq 2)$, $p>1$, $g=|x|^{\frac{2p}{N}\beta}(dx_1^2+\cdots+dx_N^2)$ be a conical metric on $\mathbb{B}$,
and $\lambda_p(\mathbb{B})=\inf\left\{\int_\mathbb{B}|\nabla u|^Ndx: u\in W_0^{1,N}(\mathbb{B}),\,\int_\mathbb{B}|u|^pdx=1\right\}$.
We prove that for any $\beta\geq 0$ and $\alpha<(1+\frac{p}{N}\beta)^{N-1+\frac{N}{p}}\lambda_p(\mathbb{B})$, there exists a constant $C$ 
such that for all radially
symmetric functions $u\in W_0^{1,N}(\mathbb{B})$ with $\int_\mathbb{B}|\nabla u|^Ndx-\alpha(\int_\mathbb{B}|u|^p|x|^{p\beta}dx)^{N/p}\leq 1$, there holds
$$\int_\mathbb{B}e^{\alpha_N(1+\frac{p}{N}\beta)|u|^{\frac{N}{N-1}}}|x|^{p\beta}dx\leq C,$$
where $|x|^{p\beta}dx=dv_g$, $\alpha_N=N\omega_{N-1}^{1/(N-1)}$, $\omega_{N-1}$ is the area of the unit sphere in $\mathbb{R}^N$; moreover, extremal functions for such inequalities exist. The case $p=N$, $-1<\beta<0$ and $\alpha=0$ was considered by Adimurthi-Sandeep \cite{A-S}, while the case $p=N=2$,
 $\beta\geq 0$ and $\alpha=0$ was studied by de Figueiredo-do \'O-dos Santos  \cite{F-do-dos}.
 \end{abstract}

\begin{keyword}
 Trudinger-Moser inequality \sep blow-up analysis \sep conical metric

\MSC[2010] 35J15; 46E35.

\end{keyword}

\end{frontmatter}

\section{Introduction}
Let $\Omega$  be a smooth bounded domain in $\mathbb{R}^N$ $(N\geq 2)$, $W_0^{1,N}(\Omega)$ be the
completion of $C_0^\infty(\Omega)$ under the Sobolev norm
$$\|u\|_{W_0^{1,N}(\Omega)}=\le(\int_\Omega|\nabla u|^Ndx\ri)^{1/N},$$
where $\nabla$ denotes the gradient operator. Write $\alpha_N=N\omega_{N-1}^{1/(N-1)}$, where $\omega_{N-1}$ stands for
 the area of the unit sphere in $\mathbb{R}^N$. Then the classical Trudinger-Moser inequality \cite{24,19,17,22,14} says
\be\label{T-M}\sup_{u\in W_0^{1,N}(\Omega),\,\|u\|_{W_0^{1,N}(\Omega)}\leq 1}\int_\Omega e^{\alpha |u|^{\f{N}{N-1}}}dx<\infty,\quad \forall
\alpha\leq \alpha_N.\ee
This inequality is sharp in the sense that if $\alpha>\alpha_N$, all integrals in (\ref{T-M}) are still finite, but the supremum is infinite.
While the existence of extremal functions for it was solved by Carleson-Chang \cite{CC}, Flucher \cite{Flucher} and Lin \cite{Lin}.

Using a symmetrization argument and a change of variables, Adimurthi-Sandeep \cite{A-S} generalized (\ref{T-M}) to the following
singular version:
\be\label{Adi-Sandeep}\sup_{u\in W_0^{1,N}(\Omega),\,\|u\|_{W_0^{1,N}(\Omega)}\leq 1}\int_\Omega \f{e^{\alpha_N\gamma|u|^{\f{N}{N-1}}}}{|x|^{N\beta}}dx
<\infty,\quad \forall\,\, 0\leq\beta<1,\,\,\,0<\gamma\leq 1-\beta.\ee
Also, it is sharp in the sense that if $\gamma>1-\beta$, integrals are still finite, but the above supremum is infinite.
The inequality (\ref{Adi-Sandeep}) was extended to the whole Euclidean space $\mathbb{R}^N$ by Adimurthi-Yang \cite{Adi-Yang}.
The existence of extremal functions for (\ref{Adi-Sandeep}) in the case $N=2$ was due to Casto-Roy \cite{C-R}, Yang-Zhu \cite{Yang-Zhu-JFA} and
Iula-Mancini \cite{I-M}. An interesting question is whether or not (\ref{Adi-Sandeep}) still holds for $\beta<0$. Generally, the answer is negative.
To see this, we choose $x_0\not= 0$, $r_0>0$ such that $B_{2r_0}(x_0)\subset \Omega\setminus\{0\}$. For any $0<\epsilon< r_0$, we write the
Moser function
$$m_\epsilon(x)=\le\{\begin{array}{lll}
\f{1}{\omega_{N-1}^{1/N}}\le(\log\f{r_0}{\epsilon}\ri)^{\f{N-1}{N}},&{\rm when} &|x-x_0|\leq \epsilon\\[1.5ex]
\f{1}{\omega_{N-1}^{1/N}}\f{\log\f{r_0}{|x-x_0|}}{\le(\log\f{r_0}{\epsilon}\ri)^{1/N}},&{\rm when}&\epsilon<|x-x_0|\leq
r_0\\[1.5ex] 0,&{\rm when} &|x-x_0|>r_0. \end{array}\ri.$$
An easy computation shows $\|m_\epsilon\|_{W_0^{1,N}(\Omega)}=1$. Since $\beta<0$, we have
$$\int_\Omega e^{\alpha_N(1-\beta)m_\epsilon^{\f{N}{N-1}}}|x|^{-N\beta}dx\geq \int_{B_\epsilon(x_0)}e^{\alpha_N(1-\beta)m_\epsilon^{\f{N}{N-1}}}|x|^{-N\beta}dx\ra \infty \quad{\rm as}\quad \epsilon\ra 0.$$
 Even worse, the above estimate still holds if $\alpha_N(1-\beta)$ is replaced by any $\alpha>\alpha_N$.
 In conclusion, the singular Trudinger-Moser inequality (\ref{Adi-Sandeep}) does not hold for $\beta<0$.

Let us consider the unit ball $\mathbb{B}\subset\mathbb{R}^N$, which is centered at the origin. Let
$\mathscr{S}$ be a set of all radially symmetric functions. With a slight abuse of notations, we say that $u$ is radially symmetric
if $u(x)=u(|x|)$ for almost every $x\in \mathbb{B}$. It was proved by Ni \cite{Ni} that $W_0^{1,N}(\mathbb{B})\cap\mathscr{S}$ can be imbedded in $L^p(\mathbb{B},|x|^\alpha)$ with $\alpha>0$ and $p=2(N+\alpha)/(N-2)$ greater than $2^\ast=2N/(N-2)$.
Motivated by results of Bonheure-Serra-Tarallo \cite{B-S-T}, Calanchi-Terraneo \cite{C-T}
 and de Figueiredo-dos Santos-Miyagaki \cite{F-dos}, de Figueiredo-do \'O-dos Santos \cite{F-do-dos} observed that
in the case $N=2$,
\be\label{supercritical}\sup_{u\in W_0^{1,2}(\mathbb{B})\cap\mathscr{S},\, \|u\|_{W_0^{1,2}(\mathbb{B})}\leq 1}
\int_{\mathbb{B}}e^{\alpha u^2}|x|^{2\gamma}dx<\infty,\quad \forall\,\alpha\leq 4\pi(1+\gamma),\,\, \gamma\geq 0;\ee
moreover, extremal function for the above supremum exists. Of course they discussed more general weight $h(|x|)$
and fast growth $F(u)$ instead of $|x|^\alpha$ and $e^{4\pi(1+\alpha)u^2}$ respectively.

Our aim is to generalize (\ref{supercritical}) to higher dimensional case and to stronger versions. We first have the following:

\begin{theorem}\label{Thm1}
Let $\mathbb{B}$ be the unit ball in $\mathbb{R}^N$ $(N\geq 2)$, $W_0^{1,N}(\mathbb{B})$ and $\mathscr{S}$ be as above. Then there holds
for any $\beta\geq 0$,
\be\label{Stand}\sup_{u\in W_0^{1,N}(\mathbb{B})\cap\mathscr{S},\,\|u\|_{W_0^{1,N}(\mathbb{B})}\leq 1}\int_\mathbb{B}
e^{\gamma|u|^{\f{N}{N-1}}}|x|^{N\beta}dx<\infty,\quad\forall\, \gamma\leq\alpha_N(1+\beta).\ee
Here $\alpha_N(1+\beta)$ is the best constant in the sense that if $\gamma>\alpha_N(1+\beta)$, all integrals are finite but
the supremum is infinity. Moreover, for any $\beta\geq 0$ and any $\gamma\leq \alpha_N(1+\beta)$,
the supremum in (\ref{Stand}) can be attained.
\end{theorem}

By a rearrangement argument, for any $\gamma\leq\alpha_N$,  there holds
\be\label{symmetry}\sup_{u\in W_0^{1,N}(\mathbb{B}),\,\|u\|_{W_0^{1,N}(\mathbb{B})}\leq 1}\int_\mathbb{B}
e^{\gamma|u|^{\f{N}{N-1}}}dx=\sup_{u\in W_0^{1,N}(\mathbb{B})\cap\mathscr{S},\,\|u\|_{W_0^{1,N}(\mathbb{B})}\leq 1}\int_\mathbb{B}
e^{\gamma|u|^{\f{N}{N-1}}}dx.\ee
Therefore, when $\Omega=\mathbb{B}$, Theorem \ref{Thm1} includes the classical Trudinger-Moser inequality (\ref{T-M}) as a special case and complements
Adimurthi-Sandeep's inequality (\ref{Adi-Sandeep}).

Motivated by \cite{Yang-JFA-06,ZhuJ,Yuan-Zhu}, we would generalize Theorem \ref{Thm1} to a version involving eigenvalue of the $N$-Laplace. For $p>1$,
define
\be\label{eigen}\lambda_p(\mathbb{B})=\inf_{u\in W_0^{1,N}(\mathbb{B}),\,u\not\equiv0}
\frac{\int_\mathbb{B}|\nabla u|^Ndx}{(\int_{\mathbb{B}}|u|^pdx)^{N/p}}.\ee
For $\alpha<\lambda_p(\mathbb{B})$, we write for simplicity
\be\label{norm}\|u\|_{1,\alpha}=\le(\int_\mathbb{B}|\nabla u|^Ndx-\alpha(\int_\mathbb{B}|u|^pdx)^{N/p}\ri)^{1/N}.\ee
\begin{theorem}\label{Thm2} Given $p>1$.
In addition to the assumptions of Theorem \ref{Thm1}, let $\lambda_p(\mathbb{B})$ and $\|\cdot\|_{1,\alpha}$ be defined as in (\ref{eigen}) and (\ref{norm})
respectively. Then if $\alpha<\lambda_p(\mathbb{B})$, there holds
$$\sup_{u\in W_0^{1,N}(\mathbb{B}),\,\|u\|_{1,\alpha}\leq 1}\int_\mathbb{B} e^{\gamma |u|^{\f{N}{N-1}}}dx<\infty,\quad
\forall\, \gamma\leq \alpha_N.$$
Moreover, the above supremum can be attained.
\end{theorem}
When $p=N$, Theorem \ref{Thm2} was proved by Nguyen \cite{Nguyen2} for a smooth bounded domain.  As a consequence of Theorem \ref{Thm2}, we improve Theorem \ref{Thm1} as follows:
\begin{theorem}\label{Thm3}
Given $p>1$.
Under the same assumptions of Theorem \ref{Thm2},  for any $\beta\geq 0$ and any $\alpha<(1+\frac{p}{N}\beta)^{N-1+N/p}\lambda_p(\mathbb{B})$ ,
there holds
\be\label{sup-2}\sup_{u\in W_0^{1,N}(\mathbb{B})\cap\mathscr{S},\,\int_\mathbb{B}|\nabla u|^Ndx-\alpha(\int_\mathbb{B}|u|^p|x|^{p\beta}dx)^{N/p}\leq 1
}\int_\mathbb{B} e^{\gamma |u|^{\f{N}{N-1}}}|x|^{p\beta}dx<\infty,\quad
\forall\, \gamma\leq \alpha_N(1+\frac{p}{N}\beta),\ee
where $\alpha_N(1+\frac{p}{N}\beta)$ is the best constant in the same sense as in Theorem \ref{Thm1}. Furthermore, the supremum in (\ref{sup-2}) can be attained.
\end{theorem}

We now explain the geometric meaning of the term $|x|^{p\beta}dx$.
Let $g_0$ be the standard Euclidean metric, namely $g_0(x)=d{x_1}^2+\cdots+d{x_N}^2$. Define a metric $g(x)=|x|^{\frac{2p}{N}\beta}g_0(x)$ for
$x\in\mathbb{B}$. Then $(\mathbb{B}, g)$ is a conical manifold with the volume element $dv_g=|x|^{p\beta}dx$. Moreover, $|\nabla u|^Ndx=|\nabla_gu|^Ndv_g$.\\

The proof of Theorems \ref{Thm1} and \ref{Thm3} is based on a change of variables. While the proof of Theorem \ref{Thm2} is based on
blow-up analysis.
In the remaining part of this note, we shall prove Theorems \ref{Thm1}-\ref{Thm3} respectively.

\section{Proof of Theorem \ref{Thm1}}
Let $\beta\geq 0$ and $\gamma\leq \alpha_N(1+\beta)$. Write for simplicity $u(x)=u(r)$ with $r=|x|$. Following
Smets-Willem-Su \cite{Su} and Adimurthi-Sandeep \cite{A-S},
we make a change of variables. Define a function
$$v(r)=(1+\beta)^{1-1/N}u(r^{1/(1+\beta)}).$$ A straightforward calculation shows
\bea\nonumber\int_\mathbb{B}|\nabla v|^Ndx&=&\omega_{N-1}\int_0^1 |v^\prime(r)|^Nr^{N-1}dr\\\nonumber
&=&\f{\omega_{N-1}}{1+\beta}\int_0^1 |u^\prime(r^{1/(1+\beta)})|^Nr^{{(N-1-\beta)/(1+\beta)}}dr\\
&=&\omega_{N-1}\int_0^1 |u^\prime(t)|^Nt^{N-1}dt=\int_\mathbb{B}|\nabla u|^Ndx\label{grad}
\eea
and
\bea\nonumber
\int_\mathbb{B}e^{\gamma|u|^{\f{N}{N-1}}}|x|^{N\beta}dx&=&\omega_{N-1}\int_0^1e^{\gamma|u(r)|^{\f{N}{N-1}}}r^{N-1+N\beta}dr\\\nonumber
&=&\omega_{N-1}\int_0^1e^{\f{\gamma}{1+\beta}|v(r^{1+\beta})|^{\f{N}{N-1}}}r^{N-1+N\beta}dr\\\label{exp}
&=&\f{\omega_{N-1}}{1+\beta}\int_0^1e^{\f{\gamma}{1+\beta}|v(t)|^{\f{N}{N-1}}}t^{N-1}dt=
\f{1}{1+\beta} \int_\mathbb{B}e^{\f{\gamma}{1+\beta}|v|^{\f{N}{N-1}}}dx.
\eea
Then it follows from (\ref{grad}), (\ref{exp}) and (\ref{symmetry}) that
\be\label{sup-equ}
\sup_{u\in W_0^{1,N}(\mathbb{B})\cap \mathscr{S},\,\|u\|_{W_0^{1,N}(\mathbb{B})}\leq 1}\int_\mathbb{B}
e^{\gamma|u|^{\f{N}{N-1}}}|x|^{N\beta}dx
=\f{1}{1+\beta}\sup_{v\in W_0^{1,N}(\mathbb{B}),\,\|v\|_{W_0^{1,N}(\mathbb{B})}\leq 1}\int_\mathbb{B}e^{\f{\gamma}{1+\beta}|v|^{\f{N}{N-1}}}dx.
\ee
According to Carleson-Chang \cite{CC}, the supremum on the right-hand side of (\ref{sup-equ}) can be attained, so does the supremum
on the left-hand side. This concludes  Theorem \ref{Thm1}. $\hfill\Box$

\section{Proof of Theorem \ref{Thm2}}

In this section, we use the standard blow-up analysis to prove Theorem \ref{Thm2}. This method was originally introduced by
Ding-Jost-Li-Wang \cite{DJLW} and Li \cite{Lijpde,Liscience},  and extensively employed by Yang \cite{Yang-JFA-06,Yang-Corrigendum,Yang-Trans,Yang-JDE-15}, Lu-Yang \cite{Lu-Yang}, Li-Ruf \cite{Li-Ruf},
Zhu \cite{ZhuJ},
do \'O-de Souza \cite{do-de-1,do-de-2},
Li-Yang \cite{Li-Yang}, Li \cite{Li}, Nguyen \cite{Nguyen,Nguyen2}
 and others. Comparing with the case $p\leq N$ \cite{ZhuJ,Nguyen2}, we need more analysis to deal with the general case $p>1$.

\subsection{The existence of maximizers for subcritical functionals}

Let $\alpha<\lambda_p(\mathbb{B})$. Denote
$$\Lambda_{\gamma,\alpha}=\sup_{u\in W_0^{1,N}(\mathbb{B}),\,\|u\|_{1,\alpha}\leq 1}\int_\mathbb{B}e^{\gamma |u|^{\f{N}{N-1}}}dx.$$
\begin{lemma}\label{subc}
For any positive integer $k$, there exists a decreasing radially symmetric function $u_k\in W_0^{1,N}(\mathbb{B})\cap C^1(\overline{\mathbb{B}})$ with $\|u_k\|_{1,\alpha}=1$
such that
$\int_\mathbb{B}e^{\gamma_k|u_k|^{\f{N}{N-1}}}dx=\Lambda_{\gamma_k,\alpha}$,
where $\gamma_k=\alpha_N-1/k$. Moreover, $u_k$ satisfies the Euler-Lagrange equation
\be\label{E-L}\le\{\begin{array}{lll}
 -\Delta_Nu_k-\alpha \left(\int_{\mathbb{B}}u_k^{p}dx\right)^{\frac{N}{p}-1}u_k^{p-1}=\f{1}{\lambda_k}u_k^{\f{1}{N-1}}e^{\gamma_ku_k^{\f{N}{N-1}}}\\[1.2ex]
 u_k\geq 0\quad{\rm in}\quad \mathbb{B}\\[1.2ex]
 u_k=0\quad{\rm on}\quad \p\mathbb{B}\\[1.2ex]
 \lambda_k=\int_\mathbb{B}u_k^{\f{N}{N-1}}e^{\gamma_k u_k^{\f{N}{N-1}}}dx,
\end{array}\ri.\ee
where $\Delta_Nu_k={\rm div} (|\nabla u_k|^{N-2}\nabla u_k)$.
\end{lemma}

\proof Let $k$ be a positive integer. By a rearrangement argument, there exists a sequence of decreasing radially symmetric functions $u_j\in W_0^{1,N}(\mathbb{B})$
with $\|u_j\|_{1,\alpha}\leq 1$ and $\int_\mathbb{B}e^{\gamma_k|u_j|^{\f{N}{N-1}}}dx\ra\Lambda_{\gamma_k,\alpha}$ as $j\ra\infty$. Since
$\alpha<\lambda_p(\mathbb{B})$, $u_j$ is bounded in $W_0^{1,N}(\mathbb{B})$. Without loss of generality, we assume
$u_j$ converges to some function $u_k$ weakly in $W_0^{1,N}(\mathbb{B})$, strongly in $L^s(\mathbb{B})$ for any $s>1$ and almost everywhere
in $\mathbb{B}$. If $u_k\equiv 0$, then $\|u_j\|_{W_0^{1,N}(\mathbb{B})}\leq 1+o_j(1)$. Thus $e^{\gamma_k u_j^{{N}/{(N-1)}}}$ is bounded in
 $L^{q}(\mathbb{B})$ for some $q>1$. It follows that $e^{\gamma_k u_j^{{N}/{(N-1)}}}$ converges to $1$ in $L^1(\mathbb{B})$.
 This implies that $\Lambda_{\gamma_k,\alpha}=|\mathbb{B}|$, the volume of $\mathbb{B}$, which is impossible. Therefore $u_k\not\equiv 0$. Clearly $u_k$ is also decreasingly radially symmetric and $\|u_k\|_{1,\alpha}\leq 1$. Define a function sequence
$$v_j=\f{u_j}{\left(1+\alpha(\int_\mathbb{B} u_j^pdx)^{N/p}\right)^{1/N}}.$$
It follows that $\|v_j\|_{W_0^{1,N}(\mathbb{B})}\leq 1$, $v_j$ converges to $v_k=u_k/(1+\alpha(\int_\mathbb{B} u_k^pdx)^{N/p})^{1/N}$
weakly in $W_0^{1,N}(\mathbb{B})$. By a result of Lions (\cite{Lions}, Theorem I.6), for any $q<1/(1-\|v_k\|_{W_0^{1,N}(\mathbb{B})}^N)^{1/(N-1)}$, there holds
\be\label{Lions-lemma}\lim_{j\ra\infty}\int_\mathbb{B}e^{q\alpha_Nv_j^{\f{N}{N-1}}}dx<\infty.\ee
One can easily check that
\be\label{strict}\le(1+\alpha(\int_\mathbb{B} u_k^pdx)^{N/p}\ri)\le(1-\|v_k\|_{W_0^{1,N}(\mathbb{B})}^N\ri)=1-\|u_k\|_{1,\alpha}^N<1.\ee
It follows from (\ref{Lions-lemma}) and (\ref{strict}) that $e^{\gamma_ku_j^{N/(N-1)}}$ is bounded in $L^r(\mathbb{B})$ for some $r>1$,
and thus $e^{\gamma_ku_j^{N/(N-1)}}\ra e^{\gamma_ku_k^{N/(N-1)}}$ in $L^1(\mathbb{B})$ as $j\ra\infty$. Hence
$\int_\mathbb{B}e^{\gamma_ku_k^{N/(N-1)}}dx=\Lambda_{\gamma_k,\alpha}$ and $u_k$ is the desired extremal function.
Clearly $\|u_k\|_{1,\alpha}=1$. Moreover, the Euler-Lagrange equation of $u_k$ is (\ref{E-L}). According to the regularity theory for degenerate elliptic equations, see Serrin \cite{Serrin}, Tolksdorf \cite{Tolksdorf} and Lieberman \cite{Lieberman}, we have
$u_k\in C^1(\overline{\mathbb{B}})$. $\hfill\Box$\\

It is indicated by Lemma \ref{subc} that for any $\gamma<\alpha_N$ and $\alpha<\lambda_p(\mathbb{B})$, the supremum $\Lambda_{\gamma,\alpha}$ can be attained.
In particular, for any $\gamma_k=\alpha_N-1/k$, there exists a maximizer $u_k\geq 0$ satisfies (\ref{E-L}). It is not difficult to see that
\be\label{limit}\lim_{k\ra\infty}\int_\mathbb{B}e^{\gamma_k u_k^{\f{N}{N-1}}}dx=\Lambda_{\alpha_N,\alpha}=
\sup_{u\in W_0^{1,N}(\mathbb{B}),\,\|u\|_{1,\alpha}\leq 1}\int_\mathbb{B}e^{\alpha_N |u|^{\f{N}{N-1}}}dx.\ee

Since $\|u_k\|_{1,\alpha}=1$, without loss of generality, we can assume that $u_k$ converges to $u_0$
weakly in $W_0^{1,N}(\mathbb{B})$, strongly in $L^s(\mathbb{B})$ for any $s>1$, and almost everywhere in $\mathbb{B}$.
Let $c_k=u_k(0)=\max_{\mathbb{B}}u_k$. If $c_k$ is bounded, then applying the Lebesgue dominated convergence theorem to
(\ref{limit}), we know that $u_0$ is the desired extremal function for the supremum $\Lambda_{\alpha_{N},\alpha}$.
Hereafter we assume \be\label{infty}c_k\ra\infty\quad {\rm as}\quad k\ra\infty.\ee

\begin{lemma}\label{concentraion}
Let $u_0$ be the limit of $u_k$ as above. Then
$u_0\equiv 0$ and $|\nabla u_k|^Ndx\rightharpoonup \delta_0$ weakly in the sense of measure, where $\delta_0$ stands for the Dirac measure
centered at the origin.
\end{lemma}

\proof We first prove $u_0\equiv 0$. Suppose not. It follows from Lions' lemma  that
$e^{\gamma_ku_k^{N/(N-1)}}$ is bounded in $L^q(\mathbb{B})$ for some $q>1$. Then applying elliptic estimates to (\ref{E-L}),
we conclude $u_k$ is uniformly bounded in $\mathbb{B}$, which contradicts our assumption (\ref{infty}). Therefore $u_0\equiv 0$.

Next we prove $|\nabla u_k|^Ndx\rightharpoonup \delta_0$. Suppose not. Since $\|u_k\|_{1,\alpha}=1$ and $u_0\equiv 0$, there would hold
$\|u_k\|_{W_0^{1,N}(\mathbb{B})}=1+o_k(1)$. Thus there exists some $0<r_0<1$ such that
$$\limsup_{k\ra\infty}\int_{|x|\leq r_0}|\nabla u_k|^Ndx<1.$$
It follows from the classical Trudinger-Moser inequality (\ref{T-M}) that $e^{\gamma_k(u_k-u_k(r_0))^{N/(N-1)}}$ is bounded in $L^q({B}_{r_0})$
for some $q>1$. Since $u_k$ is decreasing radially symmetric and $\|u_k\|_{1,\alpha}=1$ with $\alpha<\lambda_p(\mathbb{B})$, we have
$$u_k^N(r_0)\leq \left(\f{1}{|B_{r_0}|}\int_{|x|\leq r_0}u_k^pdx\right)^{N/p}\leq \f{1}{(\lambda_p(\mathbb{B})-\alpha)|B_{r_0}|^{N/p}}.$$
Hence $e^{\gamma_ku_k^{N/(N-1)}}$ is also bounded in $L^{q_1}({B}_{r_0})$
for some $q_1>1$. Then applying elliptic estimates to (\ref{E-L}), we conclude that $u_k$ is uniformly bounded in $\mathbb{B}$. This contradicts
(\ref{infty}) and ends the proof of the lemma. $\hfill\Box$

\subsection{Blow-up analysis}

Let
$r_k=\lambda_k^{\f{1}{N}}c_k^{-\f{1}{N-1}}e^{-\f{\gamma_k}{N}c_k^{{N}/{(N-1)}}}$.
Using the same argument as in the proof of (\cite{Yang-JFA-06}, Lemma 4.3), one has by Lemma \ref{concentraion} and the classical
Trudinger-Moser inequality (\ref{T-M}) that
\begin{align}\label{scale}
r_ke^{a u_k^{N/(N-1)}}\ra0~~~~\mbox{as}~~k\ra\infty,~~\forall a<\alpha_N/N.
\end{align}
For $x\in B_{r_k^{-1}}$, we define
$\psi_k(x)=c_k^{-1}u_k(r_kx)$ and $\varphi_k(x)=c_k^{\f{1}{N-1}}(u_k(r_kx)-c_k)$.
\begin{lemma}\label{lemma6}
Up to a subsequence, there holds $\psi_k\ra 1$ in $C^1_{\rm loc}(\mathbb{R}^N)$ and $\varphi_k\ra \varphi$ in
$C^1_{\rm loc}(\mathbb{R}^N)$ as $k\ra\infty$, where
\be\label{bubble-varphi}\varphi(x)=-\f{N-1}{\alpha_N}\log\le(1+\f{\alpha_N}{N^{N/(N-1)}}|x|^{\f{N}{N-1}}\ri).\ee
\end{lemma}
\proof
A simple calculation gives
\begin{align}\label{blowup}
-\Delta_N \psi_k=\alpha c_k^{p-N}r_k^N||u_k||_p^{N-p}\psi_k^{p-1}+c_k^{-N}e^{\gamma_k(u_k^{N/(N-1)}-c_k^{N/(N-1)})}\psi_k^{1/(N-1)}
\end{align}
and
\begin{align}\label{blowup2}
&-\Delta_N \varphi_k=\alpha c_k^{p}r_k^N||u_k||_p^{N-p}\psi_k^{p-1}+e^{\gamma_k(u_k^{N/(N-1)}-c_k^{N/(N-1)})}\psi_k^{1/(N-1)}.
\end{align}
Since $u_k$ is bounded in $L^p(\mathbb{B})$, one has by (\ref{infty}) and (\ref{scale}) that
\begin{align}\label{blowup-1}
\left(\int_{B_{r_k^{-1}}}\left(c_k^{p-N}r_k^N||u_k||_p^{N-p}\psi_k^{p-1}\right)^{p/(p-1)}dx\right)^{(p-1)/p}
=c_k^{1-N}r_k^{N/p}||u_k||_p^{N-1}\ra0~~\mbox{as}~~k\ra\infty.
\end{align}
Since $0\leq\psi_k\leq1$, there holds
\begin{align}\label{blowup-2}
c_k^{-N}e^{\gamma_k(u_k^{N/(N-1)}-c_k^{N/(N-1)})}\psi_k^{1/(N-1)}\ra0~~\mbox{as}~~k\ra\infty.
\end{align}
It follows from (\ref{blowup-1}) and (\ref{blowup-2}) that $\Delta_N\psi_k$ is bounded in $L^{{p}/{(p-1)}}(B_{r_k^{-1}})$.
Applying the regularity theory \cite{Tolksdorf} to (\ref{blowup}), one obtains
\begin{align}\label{psi}
\psi_k\ra\psi~~\mbox{in}~~C^0_{\rm loc}(\mathbb{R}^N)~~\mbox{as}~~k\ra\infty,
\end{align}
where
\be\label{psi-scope}\psi(0)=1,\quad 0\leq\psi(x)\leq 1,\,\,\forall x\in \mathbb{R}^N.\ee
When $1<p\leq N$, one can easily see that
\begin{align}\label{blowup2-1}
c_k^{p}r_k^N||u_k||_p^{N-p}\psi_k^{p-1}\ra 0
\end{align}
uniformly in $x\in B_{r_k^{-1}}$ as $k\ra\infty$.
When $p>N$, we have for any $R>0$ and sufficiently large $k$
\begin{align}\label{nor-eps}
||u_k||_p^{N-p}=\left(\int_{\mathbb{B}}u_k^{p}dx\right)^{N/p-1}
\leq \left(\int_{B_{Rr_k}}u_k^{p}dx\right)^{N/p-1}
=c_k^{N-p}r_k^{N^2/p-N}\left(\int_{B_R}\psi_k^pdx\right)^{N/p-1}.
\end{align}
In view of (\ref{psi-scope}), we conclude
 $\int_{B_R}\psi^pdx>0$, which together with (\ref{psi}) and (\ref{nor-eps}) leads to
 \begin{align*}
||u_k||_p^{N-p}\leq 2\le(\int_{B_R}\psi^pdx\ri)^{N/p-1}c_k^{N-p}r_k^{N^2/p-N}
\end{align*}
for sufficiently large $k$.
This together with (\ref{scale}) gives
\begin{align}\label{blowup2-2}
c_k^{p}r_k^N||u_k||_p^{N-p}\psi_k^{p-1}\leq 2\le(\int_{B_R}\psi^pdx\ri)^{N/p-1}c_k^{N}r_k^{N^2/p}
\ra0~~\mbox{as}~~k\ra\infty.
\end{align}
It then follows from (\ref{blowup2-1}) and (\ref{blowup2-2}) that $\Delta_N\psi_k$ is bounded in $L^\infty(B_R)$. Applying again the regularity theory \cite{Tolksdorf} to (\ref{blowup}),
we conclude that $\psi_k\ra \psi$ in $C^1(B_{R/2})$. Since $R$ is arbitrary, up to a subsequence, there holds
$$\label{psi-3}\psi_k\ra\psi~~\mbox{in}~~C^1_{\rm loc}(\mathbb{R}^N)~~\mbox{as}~~k\ra\infty,
$$
where $\psi$ is a solution of
\begin{align*}
-\Delta_N \psi=0~~\mbox{in}~~\mathbb{R}^N,~~~~0\leq\psi\leq\psi(0)=1.
\end{align*}
The Liouville theorem implies that $\psi\equiv1$ in $\mathbb{R}^N$.

Recalling (\ref{blowup2-1}) and (\ref{blowup2-2}), $c_k^{p}r_k^N||u_k||_p^{N-p}\psi_k^{p-1}\ra 0$ in $L^\infty_{\rm loc}(\mathbb{R}^N)$.
Then using the same argument as in (\cite{Li-Ruf}, Section 3) or (\cite{Li-Yang}, Lemma 17), we have by applying elliptic estimates to (\ref{blowup2}),
$$\label{con-bub}\varphi_k\ra\varphi \,\,{\rm in}\,\, C^1_{\rm loc}(\mathbb{R}^N)\,\,{\rm as}\,\,k\ra\infty,$$ where
$\varphi$ satisfies
$$\le\{\begin{array}{lll}
-\Delta_N\varphi=e^{\alpha_N\f{N}{N-1}\varphi}\quad{\rm in}\quad \mathbb{R}^N\\[1.2ex]
\sup_{\mathbb{R}^N}\varphi=\varphi(0)=0.
\end{array}\ri.$$
Observing that $\varphi$ is radially symmetric, one gets (\ref{bubble-varphi}) by solving the corresponding ordinary differential equation.
$\hfill\Box$\\

Lemma \ref{lemma6} describes the asymptotic behavior of $u_k$ near the blow-up point $0$. To know $u_k$'s behavior
away from $0$, by the same argument as in the proof of (\cite{Yang-JFA-06}, Lemma 4.11), we have that
$$\label{con-Gr}\le\{\begin{array}{lll} c_k^{{1}/{(N-1)}}u_k\rightharpoonup G &{\rm weakly\,\,in}&W_0^{1,q}(\mathbb{B}),
\,\,\,\forall 1<q<N\\[1.2ex]
c_k^{{1}/{(N-1)}}u_k\ra G&{\rm strongly\,\,in}&L^s(\mathbb{B}),\,\,\,\forall 1<s\leq Nq/(N-q)\\[1.2ex]
c_k^{{1}/{(N-1)}}u_k\ra G&\quad\quad{\rm in}&C^1_{\rm loc}(\overline{\mathbb{B}}\setminus\{0\}),\end{array}\ri.$$
 where $G$ is a distributional solution of
\be\label{Green-equation}-\Delta_NG-\alpha ||G||_p^{N-p}G^{p-1}=\delta_0\quad{\rm in}\quad \mathbb{B}.\ee
According to Kichenassamy-Veron \cite{K-V}, $G$ can be represented by
$$\label{Gr}G(x)=-\f{N}{\alpha_N}\log|x|+A_0+w(x),$$
where $A_0$ is a constant, $w\in C^\nu(\mathbb{B})$ for some $0<\nu<1$ and $w(0)=0$.
In view of (\ref{limit}),  we also have an analog of (\cite{Yang-JFA-06}, Proposition 5.2), namely
\be\label{up-bnd}\Lambda_{\alpha_N,\alpha}\leq|\mathbb{B}|+\f{\omega_{N-1}}{N}e^{\alpha_NA_0+\sum_{j=1}^{N-1}\f{1}{j}}.\ee
For its proof, since no new idea comes out, we omit the details but refer the readers to \cite{Yang-JFA-06} (see also \cite{ZhuJ,Li,Nguyen}).

\subsection{Test function computation}

In this subsection, we construct a sequence of functions to show that
\be\label{geq}\Lambda_{\alpha_N,\alpha}>|\mathbb{B}|+\f{\omega_{N-1}}{N}e^{\alpha_NA_0+\sum_{j=1}^{N-1}\f{1}{j}}.\ee
The contradiction between (\ref{geq}) and (\ref{up-bnd}) indicates that $c_k$ is a bounded sequence, and whence the desired extremal
function exists. This completes the proof of Theorem \ref{Thm2}.

For any positive integer $k$, we set
\be\label{test}\phi_k(x)=\le\{\begin{array}{lll}
c+\f{1}{c^{1/(N-1)}}\le(-\f{N-1}{\alpha_N}\log(1+c_N|kx|^{\f{N}{N-1}})+b\ri),&|x|<\f{\log k}{k}\\[1.2ex]
\f{G}{c^{1/(N-1)}},&\f{\log k}{k}\leq |x|\leq 1,
\end{array}\ri.\ee
where $c$ and $b$ are constants, depending only on $k$, to be determined later. To ensure $\phi_k\in W_0^{1,N}(\mathbb{B})$,
we need
$$c+\f{1}{c^{1/(N-1)}}\le(-\f{N-1}{\alpha_N}\log(1+c_N|\log k|^{\f{N}{N-1}})+b\ri)=\f{G(\f{\log k}{k})}{c^{1/(N-1)}},$$
which implies that
\be\label{c}c^{\f{N}{N-1}}=G(\f{\log k}{k})+\f{N-1}{\alpha_N}\log(1+c_N(\log k)^{\f{N}{N-1}})-b.\ee
We now calculate the energy of $\phi_k$. In view of (\ref{test}), a straightforward calculation gives
$$\label{bub-ener-test}\int_{|x|<\f{\log k}{k}}|\nabla \phi_k|^Ndx=\f{1}{c^{{N}/{(N-1)}}}\f{N-1}{\alpha_N}\le\{
-\sum_{j=1}^{N-1}\f{1}{j}+\log(1+c_N(\log k)^{\f{N}{N-1}})+O((\log k)^{-\f{N}{N-1}})\ri\}.$$
By (\ref{Green-equation}) and the divergence theorem
\bna
\int_{\f{\log k}{k}\leq |x|\leq 1}|\nabla G|^Ndx&=&\int_{\f{\log k}{k}\leq |x|\leq 1}
G(-\Delta_N G)dx+\int_{|x|=\f{\log k}{k}}G|\nabla G|^{N-1}ds\\
&=&\alpha (\int_{\f{\log k}{k}\leq |x|\leq 1}G^pdx)^{N/p}+G(\f{\log k}{k})\int_{|x|=\f{\log k}{k}}|\nabla G|^{N-1}ds\\
&=&G(\f{\log k}{k})+\alpha(\int_\mathbb{B}G^pdx)^{N/p}+O\le(\f{(\log k)^{N+p-1}}{k^N}\ri).
\ena
As a consequence
\bna
\int_\mathbb{B}|\nabla\phi_k(x)|^Ndx&=&\f{1}{c^{N/(N-1)}}\le\{-\f{N-1}{\alpha_N}\sum_{j=1}^{N-1}\f{1}{j}+
\f{N-1}{\alpha_N}\log(1+c_N(\log k)^{\f{N}{N-1}})\ri.\\
&&\quad\quad\le.+G(\f{\log k}{k})+\alpha(\int_\mathbb{B}G^pdx)^{N/p}+O((\log k)^{-\f{N}{N-1}})\ri\}.
\ena
Noting also that
$$\left(\int_\mathbb{B}\phi_k^p(x)dx\right)^{N/p}=\f{1}{c^{N/(N-1)}}\le(\int_\mathbb{B}G^pdx+O(\f{(\log k)^{2N}}{k^N})\ri)^{N/p},$$
we obtain
$$
 \|\phi_k\|_{1,\alpha}^N=\f{1}{c^{N/(N-1)}}\le(-\f{N-1}{\alpha_N}\sum_{j=1}^{N-1}\f{1}{j}+\f{N-1}{\alpha_N}\log
 (1+c_N(\log k)^{\f{N}{N-1}})+G(\f{\log k}{k})+O((\log k)^{-\f{N}{N-1}})\ri).
$$
Set $\|\phi_k\|_{1,\alpha}=1$. It then follows that
\be\label{c-n}c^{\f{N}{N-1}}=\f{N}{\alpha_N}\log k-\f{N-1}{\alpha_N}\sum_{j=1}^{N-1}\f{1}{j}+\f{N-1}{\alpha_N}\log c_N+A_0+
O((\log k)^{-\f{N}{N-1}}).\ee
This together with (\ref{c}) leads to
\be\label{b}b=\f{N-1}{\alpha_N}\sum_{j=1}^{N-1}\f{1}{j}+O((\log k)^{-\f{N}{N-1}}).\ee
When $|x|<\f{\log k}{k}$, we calculate
\bea\nonumber
\alpha_N\phi_k^{\f{N}{N-1}}(x)&=&\alpha_Nc^{\f{N}{N-1}}\le(1+\f{1}{c^{N/(N-1)}}\le(-\f{N-1}{\alpha_N}\log(1+c_N|kx|^{\f{N}{N-1}})+b\ri)\ri)^{\f{N}{N-1}}\\
\nonumber&\geq&\alpha_Nc^{\f{N}{N-1}}\le(1+\f{N}{N-1}\f{1}{c^{N/(N-1)}}\le(-\f{N-1}{\alpha_N}\log(1+c_N|kx|^{\f{N}{N-1}})+b\ri)\ri)\\
&=&\alpha_Nc^{\f{N}{N-1}}+\f{N\alpha_Nb}{N-1}-{N}\log(1+c_N|kx|^{\f{N}{N-1}}).\label{e-e}
\eea
In view of (\ref{c-n}) and (\ref{b}),
\be\label{2}\alpha_Nc^{\f{N}{N-1}}+\f{N\alpha_Nb}{N-1}=N\log k+\sum_{j=1}^{N-1}\f{1}{j}+(N-1)\log c_N+\alpha_NA_0+
O((\log k)^{-\f{N}{N-1}}).\ee
Moreover, integration by parts leads to
\bea\nonumber
\int_{|x|<\f{\log k}{k}}e^{-{N}\log(1+c_N|kx|^{\f{N}{N-1}})}dx&=&k^{-N}\int_{|y|<\log k}\f{dy}{(1+c_N|y|^{\f{N}{N-1}})^N}\\\nonumber
&=&k^{-N}\int_0^{(\log k)^{\f{N}{N-1}}}\f{N-1}{N}\f{t^{N-2}dt}{(1+c_Nt)^N}\\\label{3}
&=&k^{-N}(1+O((\log k)^{-\f{N}{N-1}})).
\eea
Combining (\ref{e-e}), (\ref{2}) and (\ref{3}), we obtain
\be\label{loc}\int_{|x|<\f{\log k}{k}}e^{\alpha_N\phi_k^{\f{N}{N-1}}(x)}dx\geq \f{\omega_{N-1}}{N}e^{\alpha_NA_0+\sum_{j=1}^{N-1}\f{1}{j}}
+O((\log k)^{-\f{N}{N-1}}).\ee
Using an inequality $e^t\geq 1+t$, we have
\be\label{out}\int_{\f{\log k}{k}\leq|x|\leq 1}e^{\alpha_N\phi_k^{\f{N}{N-1}}(x)}dx\geq |\mathbb{B}|+
\f{\alpha_N\int_\mathbb{B}{G^{\f{N}{N-1}}}dx}
{c^{{N}/{(N-1)^2}}}+O((\log k)^{-\f{N}{N-1}}).\ee
Then (\ref{geq}) follows from (\ref{loc}) and (\ref{out}) immediately.

\section{Proof of Theorem \ref{Thm3}}
As in the proof of Theorem \ref{Thm1}, we set $v(r)=(1+\frac{p}{N}\beta)^{1-1/N}u(r^{1/(1+\frac{p}{N}\beta)})$. Note that
\bna
\int_\mathbb{B}|u|^p|x|^{p\beta}dx&=&\omega_{N-1}\int_0^1|u(r)|^pr^{N-1+p\beta}dr\\
&=&\f{\omega_{N-1}}{(1+\frac{p}{N}\beta)^{p-p/N}}\int_0^1|v(r^{1+\beta})|^Nr^{N-1+N\beta}dr\\
&=&\f{\omega_{N-1}}{(1+\frac{p}{N}\beta)^{p+1-p/N}}\int_0^1|v(t)|^Nt^{N-1}dt=\f{1}{(1+\frac{p}{N}\beta)^{p+1-p/N}}\int_\mathbb{B}|v|^pdx.
\ena
Similar calculations as in (\ref{grad}) and (\ref{exp}) tell us that
\bea\nonumber
&&\sup_{u\in W_0^{1,N}(\mathbb{B})\cap\mathscr{S},\,\int_\mathbb{B}|\nabla u|^Ndx-\alpha(\int_\mathbb{B}|u|^p|x|^{p\beta}dx)^{N/p}\leq 1
}\int_\mathbb{B} e^{\gamma |u|^{\f{N}{N-1}}}|x|^{p\beta}dx\\\label{t-3}
&=&\f{1}{1+\frac{p}{N}\beta}\sup_{v\in W_0^{1,N}(\mathbb{B})\cap\mathscr{S},\,\int_\mathbb{B}|\nabla v|^Ndx-\f{\alpha}{(1+\frac{p}{N}\beta)^{N-1+N/p}}
(\int_\mathbb{B}|v|^pdx)^{N/p}\leq 1
}\int_\mathbb{B} e^{\f{\gamma}{1+\frac{p}{N}\beta} |v|^{\f{N}{N-1}}}dx.
\eea
By a rearrangement argument, we have
\bea\nonumber&&\sup_{v\in W_0^{1,N}(\mathbb{B})\cap\mathscr{S},\,\int_\mathbb{B}|\nabla v|^Ndx-\f{\alpha}{(1+\frac{p}{N}\beta)^{N-1+N/p}}
(\int_\mathbb{B}|v|^pdx)^{N/p}\leq 1
}\int_\mathbb{B} e^{\f{\gamma}{1+\frac{p}{N}\beta} |v|^{\f{N}{N-1}}}dx\\\label{sy-3}
&=&\sup_{v\in W_0^{1,N}(\mathbb{B}),\,\int_\mathbb{B}|\nabla v|^Ndx-\f{\alpha}{(1+\frac{p}{N}\beta)^{N-1+N/p}}
(\int_\mathbb{B}|v|^pdx)^{N/p}\leq 1
}\int_\mathbb{B} e^{\f{\gamma}{1+\frac{p}{N}\beta} |v|^{\f{N}{N-1}}}dx.\eea
Since $\alpha<(1+\frac{p}{N}\beta)^{N-1+N/p}\lambda_p(\mathbb{B})$ and $\gamma\leq \alpha_N(1+\frac{p}{N}\beta)$, in view of (\ref{t-3}) and (\ref{sy-3}),
Theorem \ref{Thm3} follows from Theorem \ref{Thm2} immediately. $\hfill\Box$

\bigskip

{\bf Acknowledgements}. This work is partly supported by the National Science Foundation of China (Grant Nos.
  11471014, 11401575 and 11761131002).

\end{document}